\documentclass[a4paper,10pt]{amsart}
\usepackage{amsmath}
\usepackage{mathrsfs}
\usepackage{amsfonts}
\usepackage{amssymb}
\usepackage{euscript,eufrak,verbatim}
\usepackage[dvips]{graphicx}
\usepackage[final]{epsfig}
\newtheorem{theorem}{Theorem}[section]
\newtheorem{lemma}[theorem]{Lemma}

\newtheorem{definition}[theorem]{Definition}

\DeclareMathSymbol{\varnothing}{\mathord}{AMSb}{"3F}

\begin{document}

\title{An entropic characterization of long memory stationary process}
\author{Yiming Ding and Xuyan Xiang}

\maketitle

\pagestyle{myheadings}

\markboth{ Entropic characterization of long memory process}{\sc Yiming Ding
 and Xuyan Xiang}

\begin{abstract}
Long memory or long range dependency is an important phenomenon that may arise in the analysis of time series or spatial
data.  Most of the definitions of long memory of a stationary process $X=\{X_1, X_2,\cdots,\}$ are based on the second-order properties of the process.
The excess entropy of a stationary process is the summation of redundancies which relates to the rate of convergence of the conditional entropy  $H(X_n|X_{n-1},\cdots, X_1)$ to the entropy rate. It is proved that the excess entropy is identical to the mutual information between the past and the future when the entropy $H(X_1)$ is finite.
We suggest the definition that a stationary process is long memory if the excess entropy is infinite.
Since the definition of excess entropy of a stationary process requires very weak moment condition on the distribution of the process, it can be applied to processes whose distributions without bounded second moment. A significant property of excess entropy is that it is invariant under invertible transformation, which enables us to know the excess entropy of a stationary process from the excess entropy of other process. For stationary Guassian process, the excess entropy characterization of long memory relates to popular characterization well. It is proved that the excess entropy of fractional Gaussian noise is infinite if the Hurst parameter $H \in (1/2, 1)$.

\end{abstract}

\title{}




{\bf Key words:\ }
Long memory, excess entropy,  mutual information,  stationary process, fractional Gaussian noise





\section{Introduction}
\label{Introduction}

Long memory and long  range dependence   are synonymous notions, are the
phenomenon that may arise in the analysis of time series or spatial
data, and are  very important. This importance can be judged, for example,
by the very large number of publications having one of these notions
in the title, in areas such as finance \cite{Lo2001}, econometrics \cite{Robinson2003}, internet
modeling \cite{Karagiannis2004}, hydrology \cite{Painter1998}, climate studies \cite{Varotsos2006} or linguistics \cite{Alvarez-Lacalle2006}.

A stationary process is a sequence of random variables,
whose probability law is time invariant. Long Memory phenomenon relates to the rate of decay of statistical dependence of a stationary process,
with the implication that this decays more slowly than an
exponential decay, typically a power-like decay. Some self-similar
processes may exhibit long memory, but not all processes
having long memory are self-similar. When definitions are given, they vary from author to
author (the econometric survey \cite{Guegan2005} mentions 11 different definitions). Different
definitions of long memory are used for different applications. Most of the definitions of long memory appearing in literature are based on the second-order properties of a stochastic process.
Such properties include asymptotic behavior of covariances, spectral
density, and variances of partial sums. The reasons for popularity of the
second-order properties in this context are both historical and practical:
second-order properties are relatively simple conceptually and easy
to estimate from the data. The notion of long memory is discussed from a variety of
points of view, and a comprehensive review is given by Samorodnitsky~\cite{Samorodnitsky2006}.


Concentrating too much on the correlations has a
number of drawbacks \cite{Samorodnitsky2006}. Firstly, correlations provide only very limited information about the process if
the process is "not very close" to being Gaussian.
Secondly, rate of decay of correlations may change significantly after instantaneous one-to-one transformations of the process.
Finally, what to do if the variance is infinite?

Whatever the drawbacks of using correlations to measure length of memory in
the $L^2$ case, the whole approach breaks down when the variance is infinite. Some
of the proposed ways out in specific situations included computing "correlation-
like" numbers, or using instead characteristic functions by studying the rate of
convergence to zero of the difference. This approaches have met only with limited
success.

The question then arises of whether it is possible to develop a new approach to solve the number of drawbacks.
In information theory,
the \emph{excess entropy} was developed as intuitive measure of memory stored in a stationary stochastic process, which is related to the mutual information between the infinite past and the infinite future; see Section \ref{Pre}. It has a long history and is widely employed as a measure of correlation and complexity in a
variety of fields, from ergodic theory and dynamical systems to
neuroscience and linguistics,
see Ref. \cite{Crutchfield2003} and references therein for a review.


We advocate a different approach to the problem of long memory of stationary process $X=\{X_1, X_2, \cdots \}$ by \emph{excess entropy}. We stress that a stationary process is long memory if the excess entropy of $X$ is infinite. Such a characterization of long memory stationary process admits many advantages:
\begin{enumerate}
\item The definition of excess entropy requires $H(X_1)<+\infty$ rather than the second moment condition $\mathbb{E} X_1^2< \infty$, so it can be used to detect the long memory behavior of stationary process with heavy tail distribution;

\item The excess entropy is invariant under 1-1 transformation (Theorem 2);

\item It is closely related to second moment characterization if the stationary process is Gaussian (Theorem 3).
\end{enumerate}

The paper is organized as follows. In section 2 we introduce basic concepts of information theory and show some properties of excess entropy:  the excess entropy is identical to the mutual information between the past and the future, and the excess entropy is invariant under 1-1 transformation. In section 3 and section 4, we illustrate how excess entropy relates to covariance when the stationary process is Gaussian, and show that the excess entropy is infinite for the fractional Gaussian noise with Hurst parameter $H \in (1/2, 1)$. We summarize our results in section 5.



\section{ Excess Entropy and Mutual Information}
\label{Pre}


At first,we collect some basic concepts and theorems in information theory from Chapter 8 of the book \cite{CT}.
The Shannon entropy $H(X)$ of a discrete random variable $X$ taking values $x\in S$ is defined as
\begin{equation}H(X) :=-\mathbb{E}[\log P(X)]=-\sum\nolimits_{x\in S}P(x)\log P(x),\end{equation}
where the probability that $X$ takes on the particular value $x$ is written $P(x)\equiv P(X =
x)$. The joint entropy $H(X,Y)$ of a pair of discrete random
variables $(X, Y)$ with a joint distribution $p(x, y)$ is defined as
$$H(X,Y):=-\mathbb{E} \log p(X, Y)= - \sum_x\sum_y p(x, y) \log p(x, y).$$

The conditional entropy $H(X|Y)$, which is the entropy of a $X$
conditional on the knowledge of another random variable $Y$ is
$$H(Y|X):= -\mathbb{E }\log p(Y|X)=-\sum_x\sum_y p(x, y) \log p(y|x).
$$

The differential entropy $H(X)$ of a continuous random variable $X$ with density $f(x)$ is defined as
\begin{equation}H(X) := -\int_S f(x) \log f(x)dx,\end{equation}
where $S$ is the support set of the random variable.

If $X,\ Y$ have a joint density function $f(x, y)$, the joint entropy $ H(X,Y)$ of a pair of random
variables $(X, Y)$  is defined as $$H(X,Y)=-\int f(x,y)\log f(x,y) dx dy,$$
and the conditional differential entropy $H(X|Y)$ as
$$H(X|Y)=-\int f(x,y) \log f(x|y) dx dy.$$
Since in general $f(x|y)= f(x, y)/f(y)$, we can also write
$$H(X|Y) =H(X,Y)-H(Y).$$
The following chain rule with two random variables is valid for both discrete and continuous random variables \cite{CT}
$$H(X,Y) = H(X) + H(Y|X)=H(Y)+H(X|Y).$$

Consider two random variables $X$ and $Y$ with a joint probability
density function $f(x, y)$ and marginal probability density functions $f(x)$
and $f(y)$.  The mutual information $I (X; Y)$ is the relative entropy between the joint distribution and the product distribution $f(x)f(y)$:
$$I(X;Y)=\int\int f(x,y) \log \frac{f(x,y)}{f(x)f(y)} dx dy.$$
For discrete random variables $X$ and $Y$,
$$I(X;Y)= \sum_x\sum_y p(x, y) \log \frac{p(x, y)}{p(x)p(y)}.
$$
It follows that
$$
I(X;Y)=H(X)-H(X|Y)=H(Y)-H(Y|X)=H(X)+H(Y)- H(X,Y).
$$

The joint entropy of a set $X_1,X_2, \cdots ,X_n$ of random variables with density $f(x_1, x_2, \cdots, x_n)$
is defined as
$$H_X(n):= H(X_1,X_2, \cdots,X_n) =-\int  f(x_1,\cdots, x_n) \log f(x_1,\cdots, x_n) dx_1\cdots dx_n.$$
For convenience, $H_X(n)$ is denoted by $H(n)$ if the underlying process $X$ is known, $H(0)=0$.

A stochastic process $X=\{X_i\}$ is an indexed sequence of random variables. A stochastic process is stationary if the joint probability distribution does not change when shifted in time.
The following Lemma collects some properties of the block entropy sequence $H(n)$ of a stationary process $X$.
\begin{lemma} Let $X=\{X_1,\cdots, X_n,\cdots\}$ be a stationary process such that the entropy $H(X_1)$ is finite. We have the following
\begin{enumerate}
\item The chain rule:
\begin{equation}\label{chain}H(X_1,X_2, \cdots,X_n)=\sum_{i=1}^n H(X_i|X_{i-1}\cdots X_1).\end{equation}
\item The entropy gain $h(n):=H(n)-H(n-1)=H(X_n|X_{n-1},\cdots,X_1)$ is nonnegative;
\item $h(n)$ is nonincreasing and  $\lim_{n \to \infty}h(n)=h_{\mu}$;
\item $H(n)$ is concave, $i.e.,$ for $\lambda \in (0, 1)$,
 $$H(\lambda m+(1-\lambda)n) \leq \lambda H(m)+(1-\lambda)H(n)$$ provided that $m, n $ and $\lambda m+(1-\lambda)n$ are nonnegative integers;
 \item $H(n)$ is subadditive: for all nonnegative integers $m$ and $n$,
 $$
 H(n+m) \leq H(n)+H(m);$$
 \end{enumerate}
\end{lemma}
{\bf Proof:} \begin{enumerate}
\item It follows easily from chain rule in two random variables, see \cite{CT}.

\item By the chain rule, we get the first derivative of $H(n)$ is
\begin{eqnarray*}h(n):&=&H(n)-H(n-1)\\
&=&\sum_{i=1}^n H(X_i|X_{i-1}\cdots X_1)-\sum_{i=1}^{n-1} H(X_i|X_{i-1}\cdots X_1)\\
&=&H(X_n|X_{n-1},\cdots,X_1) \ge 0.
 \end{eqnarray*}
 \item By the stationarity of $X$, we have
 \begin{eqnarray*}h(n)-h(n-1)&=&H(X_n|X_{n-1},\cdots,X_1)-H(X_{n-1}|X_{n-2},\cdots,X_1)\\
 &\leq& H(X_n|X_{n-1},\cdots,X_2)-H(X_{n-1}|X_{n-2},\cdots,X_1)\\
 &=&0.
 \end{eqnarray*}
 \item By the previous statement, the second order derivative of $H(n)$ is
 $$
 (H(n)-H(n-1))-(H(n-1)-H(n-2))=h(n)-h(n-1) \le 0,
 $$
 which indicates that $H(n)$ is concave.
 \item Since $h(n)$ is nonincreasing, by the chain rule we know that
 $$
 H(n+m)=\sum_{i=1}^{n+m}h(i)=\sum_{i=1}^{n}h(i)+\sum_{i=n+1}^{n+m}h(i)\leq H(n)+\sum_{i=1}^{m}h(i)=H(n)+H(m).
 $$

\end{enumerate}
$\hfill \Box$

The  {\it entropy rate} of a stochastic process $X=\{X_i\}$, $X_i \in R$, is defined to be
$$h_{\mu}=h_{\mu}(X) :=\lim_{n \to \infty}\frac{1}{n}H(n) =\lim_{n \to \infty} \frac{1}{n}H(X_1,X_2,\cdots,X_n)$$
if the limit exists.
For a stationary stochastic process, by Lemma 1,the above limit always exists because the sequence $H(n)$ is sub-additive.

On the other hand, by Lemma 1, $h(n)$ is nonnegative and nonincreasing, which will approach to the same limit $h_{\mu}$ because
$$
h_{\mu}=\lim_{n\to \infty}\frac{H(n)}{n}=\lim_{n\to \infty}\frac{H(n)-H(n-1)}{n-(n-1)}=\lim_{n\to \infty} h(n).
$$
The existence of entropy rate for stationary process is proved by Shannon \cite{}. A significant property of entropy rate is the AEP (asymptotic
equipartition property), also known as the the Shannon-McMillan-Breiman theorem: If
$h_{\mu}$ is the entropy rate of a finite-valued stationary ergodic process $\{X_n\}$,
then
$$ \frac{1}{n}\log p(X_0,\cdots, X_{n}) \to  h_{\mu}$$
 with probability 1.
The entropy rate $h_{\mu}$ quantifies the irreducible
randomness in sequences produced by a source: the randomness
that remains after the correlations and structures in longer and longer sequence blocks are taken into account. It is known as the metric entropy in ergodic theory.

For stationary process $X$, by the definition of entropy rate,
$$H(n)\sim nh_{\mu} \ \ \ \  as \ n \to \infty.$$
However, knowing the value $h_{\mu}$ indicates nothing about how
$H(n)/n$ approaches this limit. Moreover, there may be sublinear terms in $H(n)$. For example,
one may have $H(n)\sim n h_{\mu} +c$  or$H(n)\sim nh_{\mu}+\log n$. The sublinear terms in $H(n)$
and the manner in which $H(n)$ converges to its asymptotic
form may reveal important structural properties about a process.

The {\it excess entropy} of a stationary process $X$ is defined as:
$$E=\sum_{n=1}^{\infty} (h(n)-h_{\mu}).$$
$h(n)-h_{\mu}$ is referred as per-symbol redundancy $r(n)$, because it tells us how much additional information must be
gained about the process in order to reveal the actual per-symbol
uncertainty $h_{\mu}$. In other words, the excess entropy $E$ is the summation of per-symbol redundancy \cite{Crutchfield2003}.
Substitute $h(n)=(H(n)-H(n-1))$ into the definition of the excess entropy, we  know that
$$
E=\lim_{n \to \infty} (H(n)-nh_{\mu}).
$$
If the excess entropy is finite, we obtain
$$H(n)\approx n h_{\mu} +E$$
as $n \to \infty$.

Observe that $H(n)/n \to h_{\mu}$, one may ask if $F=\sum_{n=1}^{\infty}(\frac{H(n)}{n}-h_{\mu})$ can be used to describe the long memory of stationary process like excess entropy. It is not the case. In fact, 
\begin{eqnarray*}
F&=&\sum_{n=1}^{\infty}(\frac{H(n)}{n}-h_{\mu})\\
&=&\sum_{n=1}^{\infty}\frac{(H(n)-H(n-1))+\cdots+(H(1)-H(0))-nh_{\mu}}{n}\\
&=& \sum_{n=1}^{\infty}\frac{(h(1)-h_{\mu})+(h(2)-h_{\mu})+\cdots+(h(n)-h_{\mu})}{n} \\
&\geq  & (h(1)-h_{\mu})\sum_{n=1}^{\infty}\frac{1}{n},
\end{eqnarray*}
which is divergent for all stationary process with $H(1)>h_{\mu}$. Notice that $H(1)=h_{\mu}$ if and only if the stationary process is independent, so $F=\infty$ unless the process is {\it i.i.d.}.

The mutual information $I_{p-f}(n):=I(X_1 \cdots X_n; X_{n+1} \cdots X_{2n})$ is the information between the history and future with length $n$, and the mutual information between past and future is the limit $I_{p-f}=\lim_{n \to \infty}I_{p-f}(n)$.

\begin{theorem} Let $X=\{X_1,\cdots, X_n,\cdots\}$ be a stationary process such that the entropy $H(X_1)$ is finite. Then the excess entropy and  the mutual information are identical:
$$E=I_{p-f}.$$
\end{theorem}

{\bf Proof:}
Denote $$E_n:=\sum_{k=1}^n (H(k)-H(k-1)-h_{\mu})=\sum_{k=1}^n (h(k)-h_{\mu})=H(n)-nh_{\mu}.$$
\begin{eqnarray*}I_{p-f}(n):&=&I(X_1,\cdots, X_n;X_{n+1},\cdots, X_{2n})\\
&=&H(X_1,\cdots, X_n)+H(X_{n+1},\cdots, X_{2n})-H(X_{1},\cdots, X_{2n})\\
&=&2H(n)-H(2n).\end{eqnarray*}
Put \begin{eqnarray*}D_n:&=&H(n)-nh(n)=nH(n-1)-(n-1)H(n)\\
&=&n(n-1)(\frac{H(n-1)}{n-1}-\frac{H(n)}{n}) \ge 0.\end{eqnarray*}

Denote $a_n=h(n)-h_{\mu}$, by the definition of entropy rate, we know that $a_n \ge 0$, $a_n \ge a_{n+1}$ for positive integer $n$, and $\lim_{n \to \infty}a_n=0$.
Furthermore, we have
\begin{equation}\label{eq}
E_n=\sum_{k=1}^n a_k, \ \ \ \ \ I_{p-f}(n)=\sum_{k=1}^n a_k-\sum_{k=n+1}^{2n} a_k, \ \ \ \ \ D_n=\sum_{k=1}^n a_k-na_n.
\end{equation}
We conclude that
\begin{equation}\label{ineq}
E_n\ge I_{p-f}(n)\ge D_n.
\end{equation}
In fact, the first inequality follows from $\{a_n\}$ is nonnegative, and the second inequality follows from $\{a_n\}$ is nonincreasing.

Since $E_n$ is the partial summation of the nonnegative series $\{a_n\}$, $E_n$ is nondecreasing.
Notice that
\begin{eqnarray*}I_{p-f}(n+1)-I_{p-f}(n):&=&(\sum_{k=1}^{n+1} a_k-\sum_{k=n+2}^{2n+2} a_k)-(\sum_{k=1}^n a_k-\sum_{k=n+1}^{2n} a_k)\\
&=&2a_{n+1}-(a_{2n+1}+a_{2n+2}) \\
&=&(a_{n+1}-a_{2n+1})+(a_{n+1}-a_{2n+2}) \ge 0,
\end{eqnarray*}
and
\begin{eqnarray*}D(n+1)-D_n:&=&(\sum_{k=1}^{n+1} a_k-(n+1)a_{n+1})-(\sum_{k=1}^n a_k-na_{n})\\
&=&n(a_n-a_{n+1}) \ge 0,
\end{eqnarray*}
both $I_{p-f}(n)$ and $D_n$ are nondecreasing. It follows that the following three limits exists:
\begin{equation}\label{lim}
\lim_{n \to \infty}E_n=E, \ \ \ \ \ \lim_{n \to \infty} I_{p-f}(n)=I_{p-f}, \ \ \ \ \ \lim_{n \to \infty}D_n=D.
\end{equation}

{\bf Claim:} $E_n$ is convergent if and only if $D_n$ is convergent, and they have the same limit when they are convergent.

If $E_n=\sum_{k=1}^n a_k$ is convergent, then the series $\sum_{k=1} 2^ka_{2^k}$ is convergent. So $2^ka_{2^k} \to 0$  as $k \to \infty$. Observe that for $2^k<n \leq 2^{k+1}$,$n a_n \leq 2^{k+1}a_{2^k}$, we get $\lim_{n \to \infty}na_n=0$.
Hence $$\lim_{n \to \infty}D_n=\lim_{n \to \infty} (\sum_{k=1}^n a_k-na_n)=\lim_{n \to \infty}E_n-\lim_{n \to \infty}na_n=\lim_{n \to \infty}E_n.$$

Suppose $D_n$ is convergent. Notice that
$$D_n=\sum_{k=1}^n a_k-na_n=\sum_{k=1}^{n-1}k(a_k-a_{k+1}),$$
it follows that $\forall\varepsilon>0$, $\exists N>0$, $\forall n, s \in \mathrm{N}$, $n>N$,
$$\sum_{k=n}^{n+s}k(a_k-a_{k+1}) \le \frac{\varepsilon}{2}.
$$
As a result,
$$
0 \le n(a_n-a_{n+s+1})=\sum_{k=n}^{n+s}n(a_k-a_{k+1}) \le \sum_{k=n}^{n+s}k(a_k-a_{k+1}) \le \frac{\varepsilon}{2}.
$$
Hence, $$ na_n \le \frac{\varepsilon}{2}+n a_{n+s+1}.$$
On the other hand, $\lim_{s \to \infty}a_{n+s+1}=0$ implies that $\forall n>N$, $\exists S \in \mathrm{N}$, $\forall s>S$,
$$
0<a_{n+s+1}<\frac{\varepsilon}{2n}.$$
We obtain for $n>N$,
$$
0 \le na_n \le \frac{\varepsilon}{2}+n a_{n+s+1}\le \frac{\varepsilon}{2}+\frac{\varepsilon}{2}=\varepsilon,
$$
which indicates $\lim_{n \to \infty}na_n=0$.

Therefore,
$$\lim_{n \to \infty}E_n=\lim_{n \to \infty} \sum_{k=1}^n a_k=\lim_{n \to \infty}D_n+\lim_{n \to \infty}na_n=\lim_{n \to \infty}D_n.$$
The claim is proved.

By (\ref{ineq}) and the claim, we conclude that
$$E=I_{p-f}=D.$$
$\hfill \Box$

{\bf Remarks:}
\begin{enumerate}
\item The equality $E=I_{p-f}$ is claimed in \cite{Crutchfield2003} for discrete stationary process, an heuristic "proof" was also given. The proof is too simple to be complete.

  \item The definition of excess entropy $E$ is depend on the entropy rate $h_{\mu}$. The equation $E=I_{p-f}=D$ provide two series to approximate the excess entropy $\{2H(n)-H(2n)\}_n$ and $ \{nH(n-1)-(n-1)H(n)\}_n$, which enable us to obtain lower bound of the excess entropy $E$ without knowing the entropy rate $h_{\mu}$.
\end{enumerate}

In what follows we show that the excess entropy is invariant under 1-1 transformation.

\begin{lemma}Let $Z$ be a random variable with finite entropy, $T$ is a deterministic transformation, then $H(Z) \geq H(T(Z))$. Moreover, if $T$ is an invertible transformation, we have $H(Z)=H(T(Z))$.

\end{lemma}

{\bf Proof.}
 By the chain rule of two random variables, we have
 $H(Z, T(Z))=H(Z)+H(T(Z)|Z)=H(T(Z))+H(Z|T(Z))$. Note that $H(T(Z)|Z)=0$ and $H(Z|T(Z)) \ge 0$, we know that $H(Z) \ge H(T(Z))$. Furthermore, if $T$ is invertible, we have
 $$
 H(Z) \ge H(T(Z)) \ge H(T^{-1}(T(Z)))=H(Z),
 $$
 which implies $H(Z)=H(T(Z))$.
$\hfill \Box $

\begin{theorem} Let $X=\{X_1, X_2,\cdots \}$ be a stationary stochastic process, $T$ be an invertible transformation, then the excess entropy of the process $Y:=T(X)=\{T(X_1), T(X_2), \cdots\}$  equals to  the excess entropy of $X$.
\end{theorem}

{\bf Proof.} Put $Z_n=(X_1,\cdots X_n)$, $T(Z_n)=(T(X_1),\cdots,T(X_n))$. Since $T$ is an invertible transformation, by lemma 2, we have $H(Z_n)=H(T(Z_n))$, which implies $H_X(n)=H_Y(n))$ for all $n \ge 0$.
As a result, the stationary processes $X$ and $Y$ admit the same the excess entropy.
$\hfill \Box$

To characterize the long memory by excess entropy, a new definition is suggested as follows.
\begin{definition}
\label{EEC}
A stationary process is long memory if its excess entropy is infinite.
\end{definition}

As examples, the excess entropy for two simple short memory processes are as follows.

For {\it i.i.d.} case, the entropy rate $h_{\mu}=H(1)$ because  $H(n)=n H(1)$, the excess entropy $E=0$.

For a Markov chain $X$ defined on a countable number of states, given the transition matrix $P_{ij}$, the entropy rate of $X$ is given by:
    $$\displaystyle h_{\mu}(X) = - \sum_{ij} \mu_i P_{ij} \log P_{ij}$$
where $\mu_i$ is the stationary distribution of the chain. By the Markovian property, the excess entropy of the Markov chain is:
$$E(X)=(H(1)-h_{\mu})+(H(2)-H(1)-h_{\mu})=-\sum_{i} u_i \log u_i + \sum_{ij} \mu_i P_{ij} \log P_{ij}.$$

\section{Stationary Gaussian process}

For a stationary Gaussian stochastic process, we have
$$H(X_1,X_2,\cdots, X_n) = \frac{1}{2}\log(2\pi e)^n|K^{(n)}|,$$
where the covariance matrix $K^{(n)}$ is Toeplitz with entries $r(0),r(1),\cdots,
r(n-1)$ along the top row, i.e.,  $K^{(n)}_{ij}= r(i-j) = E(X_i-EX_i)(X_j-EX_j)$.
It is useful to consider a stochastic process in the frequency domain. 
A stationary zero mean Gaussian random process is completely described by its mean correlation function $r_{k,j} = r_{k-j}$ or,
equivalently, by its power spectral density function $f$, the Fourier transform of the covariance function:
$$f(\lambda) =\sum_{n=-\infty}^{\infty}r_n \exp (in\lambda),$$
$$\gamma_k =\frac{1}{2 \pi}\int_0^{2 \pi} f(\lambda) \exp ( -i \lambda k) d\lambda.$$

Indeed, Kolmogorov showed that the entropy rate of a stationary Gaussian
stochastic process can be expressed as
$$h_{\mu}(X)=\frac{1}{2}\log (2\pi e)+\frac{1}{4\pi}\int_{-\pi}^{\pi} \log f(\lambda) d\lambda,
$$
where $S(\lambda)$ is the power spectral density of the stationary Gaussian process $X$.
On the other hand, $$h_{\mu}(X)=\lim_{n \to \infty}(H(n)-H(n-1))=\frac{1}{2}\lim_{n\to \infty}\log(2\pi e)\frac{|K^{(n)}|}{|K^{(n-1)}|}.$$

Let $X$ be a Guassian stationary process with spectral density $f(\lambda)$. Put
$$
b_k=\frac{1}{2 \pi} \int_{- \pi}^{\pi} \log f(\lambda) e^{-k \lambda} d\lambda,
$$
for any integer $k$ if it is well defined, and we refer to them as cepstrum coefficients~\cite{Lei Li 2006}.

Li proved the following theorem (Theorem 1 in~\cite{Lei Li 2006})
\begin{theorem}\label{Li-lei} Let $\{X_n\}$ be a Guassian stationary process~\cite{Lei Li 2006}\\
(i) The mutual information between the past and the future $I_{p-f}$ is finite if and only if the cepstrum coefficients satisfy the condition: $\sum_{k=-\infty}^{\infty}k b_k^2 < \infty$. In this case, $I_{p-f}=\frac{1 }{2}\sum_{k=-\infty}^{\infty}k b_k^2.$\\
(ii) If the spectral density $f(\lambda)$ is continuous, and $f(\lambda)>0$, then $I_{p-f}$ is finite if and only if the auotocovariance functions satisfy the condition $\sum_{k=-\infty}^{\infty}k \gamma_k^2<\infty.$
\end{theorem}

\begin{lemma}\label{momoton sequence}
Let $\{c(n)\}$ be a decreasing positive series. Then
$$
\sum_{n}
c(n) < \infty \ \ \ \ \ implies \ \ \  \ \ \sum_{n}
n c^2(n)<\infty.$$
\end{lemma}

{\bf Proof:} Since $\{c(n)\}$ be a decreasing positive series, $\sum_{n}
c(n)$ is convergent implies that $\{n c(n)\} \to 0$ . So $\{n c(n)\} $ is bounded by a constant $M$. Therefore,
$$
\sum_{n} n c^2(n) \leq \sum_n M c(n) =M \sum_n c(n)<\infty.
$$
$\hfill \Box$

In general, we have the following result for Gaussian stationary process.
\begin{theorem}
\label{implication}
Let $X$ be a Gaussian stationary process with decreasing autocorrelation $\{|r(n)|\}$. Suppose the spectral density $f(\lambda)$ is continuous and $f(\lambda)>0$. Then
$$
\sum_{n}
|r(n)| < \infty \Rightarrow \ \ \ the \ excess\ \ entropy\  E<\infty.$$
In other words, $X$ is not long memory in the sense of covariance implies that it is also not long memory in the sense of excess entropy.
\end{theorem}

The proof of Theorem \ref{implication} follows immediately from Theorem \ref{Li-lei} and Lemma \ref{momoton sequence}.

\noindent{\bf Remark: }
Put $a(n)=\frac{1}{n \log n}$, then $\sum_n a(n)=\infty$, but $\sum_n n a^2(n)< \infty$.
The converse implication in Lemma \ref{momoton sequence} is not true.

So for process considered in Theorem \ref{implication}, \emph{the long memory in the sense of excess entropy is more strict than that in the sense of covariance}.

\section{Fractional Gaussian Noise}
\label{FGN}

Fractional Brownian motion (FBM) with Hurst parameter $H \in (0, 1)$ is a centered continuous-time Guassian process $B^H(.)$ with covariance function
\begin{equation}
\label{ACF-1}
\rho(s,t)\equiv\mathbb{E}[B^H_tB^H_s]=\frac{1}{2}(|t|^{2H}+|s|^{2H}-|t-s|^{2H}).
\end{equation}
for $s,t\geq 0$. $B^H$ reduces to an ordinary Brownian motion for $H=1/2$. The incremental process of a FBM is a stationary discrete-time process and is called \emph{fractional Gaussian noise}, or FGN. We define of the FGN $X=\{X_k: k=0, 1,...\}$ by the autocovariance function
$$
\rho_k\equiv \mathbb{E} X_{n+k}X_n=\frac{1}{2}[|k-1|^{2H}-2|k|^{2H}+|k+1|^{2H}]
$$
for $k \in \mathbb{Z}$. It is easy to see that
\begin{equation}
\label{ACF-2}
r_k\sim H(2H-1)|k|^{-2(1-H)}.
\end{equation}
as $|k|\to\infty$. Of course, if $H=1/2$, then $\rho_k=0$ for all $k\geqslant 1$ (a Brownian motion has independent increments).
One can conclude that the summability of correlations ($\sum^{\infty}_{k=1}|\rho_k|<+\infty$) holds if $0< H<1/2$ and it does not hold if $1/2 <H<1$.
Therefore, a FGN with $H>1/2$ has become commonly accepted as having long memory, and lack of summability of correlations as a popular definition of long memory 
(Eq.
(\ref{ACF-2})).

The stationarity of the increments of the FBM implies that this is a stationary Gaussian process.

In particularly, we have the following result for FGN.
\begin{theorem}
\label{theo-FBM}
  The excess entropy of the (discrete) increment process of a fractional Brownian motion with Hurst parameter $H \in (1/2, 1)$ is infinite, i.e., it is long memory in the sense of excess entropy.
\end{theorem}




\noindent{\bf Proof:} The spectral density of the increment of fractional Brownian motion $B^H(t)$ (fractional Guassian noise) is obtained by Sinai \cite{Sinai1976}
$$
f(\lambda)=2\sin(\pi H)\Gamma(2H+1)(1-\cos \lambda)[|\lambda|^{-2H-1}+A_H(\lambda)],
$$
where $\Gamma(.)$ denotes the Gamma function and
$$
A_H(\lambda):=\sum_{j=1}^{\infty}\{(2\pi j+\lambda)^{-2H-1}+(2\pi j-\lambda)^{-2H-1}\},
$$
for\ $-\pi\leq \lambda \leq \pi$.

The spectral density $f(\lambda)=C(1-\cos \lambda)$ $\sum_{j=-\infty}^{\infty}|2\pi j+\lambda|^{-(2H+1)}$ is positive, but is not continuous at $\lambda=0$ when $H \in (1/2 ,1)$ because it is proportional to $|\lambda|^{1-2H}$ near $\lambda=0$.
Since $(1-\cos \lambda)=\sum_{n\ge 1}(-1)^{n-1}\frac{\lambda^{2n}}{(2n)!}$,
\begin{eqnarray*}
|\lambda|^{-(2H+1)}(1-\cos \lambda) &=&|\lambda|^{-(2H-1)}\sum_{n\ge 1}(-1)^{n-1}\frac{\lambda^{2(n-1)}}{(2n)!}\\
&:=&|\lambda|^{-(2H-1)}(1+g_1(\lambda))/2.
\end{eqnarray*}

Now we estimate the excess entropy of FGN via the logarithm of the spectral density.
We have
\begin{eqnarray*}
\log f(\lambda)&=&\log C + \log (1-\cos \lambda)|\lambda|^{-(2H+1)} \\
&&+\log (1+|\lambda|^{(2H+1) }A_H(\lambda) ).
\end{eqnarray*}
$\log f(\lambda)$ is an even function on $[-\pi, \pi]$, i.e., $\log f(\lambda)=\log f(-\lambda)$.
 For $n \neq 0$, we obtain the following decomposition
\begin{eqnarray*}
\gamma(n)&=&\frac{1}{2 \pi} \int_{- \pi}^{\pi} \log f(\lambda) \exp(in\lambda) d \lambda \\
&=&\frac{1}{2 \pi} \int_{- \pi}^{\pi} \log f(\lambda) \cos (n \lambda) d \lambda\\
&=& \frac{1}{2 \pi} \int_{- \pi}^{\pi} \cos (n \lambda) \log [(1-\cos \lambda)|\lambda|^{-(2H+1)}] d \lambda \\
& &+\frac{1}{2 \pi} \int_{- \pi}^{\pi} \cos (n \lambda)\log (1+  |\lambda|^{-(2H+1)} A_H(\lambda) ) d \lambda\\
&=& \frac{1}{2 \pi} \int_{- \pi}^{\pi} \cos (n \lambda) \log [|\lambda|^{-(2H-1)}(1+g_1(\lambda))/2] d \lambda \\
& & +\frac{1}{2 \pi} \int_{- \pi}^{\pi} \cos (n \lambda)\log (1+  |\lambda|^{-(2H+1)} A_H(\lambda) ) d \lambda\\
&=& \frac{1}{2 \pi} \int_{- \pi}^{\pi} \cos (n \lambda) \log |\lambda|^{-(2H-1)} d \lambda \\
& &+\frac{1}{2 \pi} \int_{- \pi}^{\pi} \cos (n \lambda) \log [(1+g_1(\lambda))/2] d \lambda\\
& & +\frac{1}{2 \pi} \int_{- \pi}^{\pi} \cos (n \lambda)\log (1+  |\lambda|^{-(2H+1)} A_H(\lambda) ) d \lambda\\
&:=& \gamma_1(n)+\gamma_2(n)+\gamma_3(n).
\end{eqnarray*}

For $\gamma_1(n)$, we get

\begin{eqnarray*}
\gamma_1(n)&=&\frac{1}{2 \pi} \int_{- \pi}^{\pi} \cos (n \lambda)\log |\lambda|^{-(2H-1)} d \lambda \\
&=& \frac{-(2H-1)}{ \pi} \int_{0}^{\pi} \cos (n \lambda)\log \lambda d \lambda \\
&=&\frac{(2H-1)}{ \pi n} \int_{0}^{\pi} \log \lambda  d \sin (n \lambda)\\
&=&\frac{-(2H-1)}{ \pi n} \int_{0}^{\pi} \frac{ \sin (n \lambda)}{\lambda} d \lambda.
\end{eqnarray*}
Observe that
\begin{equation}
\int_0^{\infty} \frac{\sin(ax)}{x}d x=\left \{ \begin{array}{ll}
\frac{\pi}{2}  & a>0 \\
-\frac{\pi}{2} & a<0.
\end{array}
\right.
\end{equation}
So $|\gamma_1(n)| \to \frac{2H-1}{2 n}$ as $n$ goes to infinity.

For the twice continuously differentiable  function $g$, the Fourier coefficient of order $n$ behaves like $O(\frac{1}{n^2})$ \cite{Stein and Shakarchi2003}.
Observe that $\log [(1+g_1(\lambda))/2]$ is twice differentiable, there exists positive constant $M_1<\infty$
$$
\gamma_2(n)=\frac{1}{2 \pi} \int_{- \pi}^{\pi} \cos (n \lambda) \log [(1+g_1(\lambda))/2] d \lambda \leq M_1/n^2.
$$

Now we estimate $\gamma_3(n)$. Denote
$$
g_2(x)=\log (1+|x|^{2H+1}A_H(x)).
$$

We have
$$
g_2'(x)=\frac{(|x|^{2H+1}A_H(x))'}{1+|x|^{2H+1}A_H(x)},
$$
{\small
$$
g_2''(x)=\frac{(|x|^{2H+1}A_H(x))''|x|^{2H+1}A_H(x)-((|x|^{2H+1}A_H(x))')^2}{(1+|x|^{2H+1}A_H(x))^2},
$$
}
where
\begin{eqnarray*}
(|x|^{2H+1}A_H(x))'&=&(2H+1)|x|^{2H}A_H(x)+|x|^{2H+1}A_H'(x) \\
(|x|^{2H+1}A_H(x))''&=&2H(2H+1)|x|^{2H-1}A_H(x)\\
& &+(2H+1)|x|^{2H}A_H'(x) \\
& &+|x|^{2H+1}A_H''(x)\\
\end{eqnarray*}
\begin{eqnarray*}
A_H'(x)&=&\sum_{j=1}^{\infty} (\frac{-(2H+1)}{(2\pi j +x)^{2H+2}}+\frac{(2H+1)}{(2\pi j -x)^{2H+2}})\\
A_H''(x)&=&\sum_{j=1}^{\infty} (\frac{(2H+1)(2H+2)}{(2\pi j +x)^{2H+3}}+\frac{(2H+1)(2H+2)}{(2\pi j -x)^{2H+3}}) .
\end{eqnarray*}

Since $|x|^{2H-1}$,$|x|^{2H}$,$|x|^{2H+1}$,$ A_H(x), A_H'(x), A_H''(x)$ are continuous  functions on
$[-\pi, \pi]$, $g_2''(x)$ are also continuous functions on $[-\pi, \pi]$. We conclude that  $$|\gamma_3(n)|\leq \frac{M}{n^2}$$ for some positive constant $M_2<\infty$.

Combine the estimations of $\gamma_1(n), \ \gamma_2(n), \ \gamma_3(n)$,
$$
\lim_{n \to \infty}\frac{|\gamma(n)|}{n}=\lim_{n\to \infty}\frac{|\gamma_(n)+\gamma_2(n)+\gamma_3(n)|}{n}=(2H+1)/2.
$$

It follows that there exists positive integer $N_0$ such that $|\gamma (n)| \geq \frac{(2H+1)}{4n}$ for $n>N_0$.

Hence, by Theorem 2, the excess entropy
\begin{eqnarray*}
E&=&\sum_{n=1}^{\infty}n |\gamma (n)|^2 >\frac{2H+1}{4} \sum_{n>N_0} n |\gamma (n)|^2 \\
&\ge&  \frac{2H+1}{4}\sum_{n>N_0} n \frac{1}{n^2}=\infty.
\end{eqnarray*}
$\hfill \Box$
\section{Conclusion}

The finiteness or infiniteness of excess entropy can be regarded as a sign between the short memory and long memory stationary processes. The definition of excess entropy of a stationary process requires very weak moment condition on the distribution of the process, and can be applied to processes whose distributions without bounded second moment. The most significant property of excess entropy is that it is invariant under invertible transformation. The invariance under invertible transformation enables us to know the excess entropy of a stationary process from the excess entropy of other process. Since conditional entropy can capture the dependence between random variables well, excess entropy is relevant to capture the dependence of a stationary process whose distribution far from Gaussian distribution. For stationary Guassian process, the excess entropy characterization of long memory related to popular characterization neatly. The challenge is to develop pertinent methods and algorithms to estimate or approximate the excess entropy of typical stationary process or  sequential data.

\section*{Acknowledgments}

We would like to thank Professor Yimin Xiao and Professor Yaozhong Hu for useful discussions and suggestions.
This work was supported by National Natural Science
Foundation of China (grant no. 11171101), and the Project Sponsored by the Scientific Research Foundation for the Returned Overseas Chinese Scholars, State Education Ministry (2015).





\vspace{0.5cm}

Yiming Ding \ \ \  ding@wipm.ac.cn
\vspace{0.5cm}

Xuyan Xiang \ \ \  xyxiang2001@126.com
\end{document}